\documentclass[preprint,11pt]{elsarticle}

\usepackage{amssymb}
\usepackage{enumerate}
\usepackage{amsmath}
\usepackage{amsfonts}
\usepackage{mathrsfs}

\usepackage[dvipsnames,usenames]{xcolor}
\usepackage{hyperref}

\hypersetup{
colorlinks = true,
urlcolor = red,
linkcolor = blue,
}

\def\R{\mathbb{R}}

\def\R{\mathbb{R}}

\def\R{\mathbb{R}}

\def\ds{\displaystyle}

\newtheorem{theorem}{Theorem}[section]

\newtheorem{e-proposition}[theorem]{Proposition}
\newtheorem{corollary}[theorem]{Corollary}
\newtheorem{e-definition}[theorem]{Definition}
\newtheorem{remark}{\it Remark\/}

\numberwithin{equation}{section}

\begin{document}

\begin{frontmatter}

 \author{Ch\'erif Amrouche $^a$}
 \ead{cherif.amrouche@univ-pau.fr}
 \address[authorlabel1]{Laboratoire de Math\'ematiques et Leurs Applications, UMR CNRS 5142\\ Universit\'e de Pau et des Pays de l'Adour,  64000 Pau,  France}
 \author{Mohand Moussaoui $^b$}
 \ead{mmohand47@gmail.com}
 \address[authorlabel2]{Laboratoire des
{\'Equations} aux D\'eriv\'ees Partielles Non Lin\'eaires et
Histoire des Math\'ematiques, Ecole Normale Sup\'erieure, Kouba, Algeria
}

\title{The Dirichlet problem for the Laplacian \\ in Lipschitz domains. Abstract}

\begin{abstract}
The main purpose of this paper is to address some questions concerning boundary value problems related to the Laplacian
 and bi-Laplacian operators, set in the framework of classical $H^s$ Sobolev spaces on a bounded Lipschitz domain of $\R^N.$ These questions are not new and a lot of work has been done in this direction by many authors using various techniques since the 80's. 
If for regular domains almost every thing is elucidated, it is not the case for Lipschitz ones and for $s$ of the form
 $s = k + 1/2$, with $k$ integer. It is well known that this framework is delicate.
Even in these cases many results are well established but sometimes not satisfactory.
Several questions remain posed.
Our main goal through this work is on one hand to give some improvements to the theory and on another one by using techniques which do not require too intricate calculations. We also tried to obtain maximal regularity for the solutions and as far as we can optimality
of the results.
\end{abstract}


\end{frontmatter}



The purpose of this work is to study the Dirichlet  problem:
$$
(\mathscr{L}_D)\ \ \ \  -\Delta u = f\quad \ \mbox{in}\ \Omega \quad
\mbox{and } \quad u = g \ \ \mbox{on }\Gamma,
$$
with data in some Sobolev spaces and the domain $\Omega$ is only Lipschitz or sometimes of class $\mathscr{C}^{1, 1}$. When $g = 0$ 
we denote this problem by $(\mathscr{L}_D^0)$ 
and when $f = 0$ we denote it by $(\mathscr{L}_D^H)$. 
These issues have been widely studied since the 1960's. In \cite{Lions}, Lions and Magenes made a complete study for smooth domains and $L^2$-theory. Grisvard in \cite{Gri} and Ne\v{c}as  \cite{Necas} treated the case where $\Omega$ is of class $\mathscr{C}^{r,1}$, with nonnegative integer $r$.  Grisvard \cite{Gri1}, \cite{Gri2} was also interested in the particular case of polygons or polyhedra. We recall that one consequence of Calder\'on-Zygmund's theory of singular integrals and boundary layer potential is that, for every $f\in W^{m-2,p}(\Omega)$ and $g\in W^{m-1/p,p}(\Gamma)$ with a positive integer $m$, the problem $(\mathscr{L}_D)$ has a unique solution $u\in W^{m,p}(\Omega)$ when $\Omega$ is of class $\mathscr{C}^{r,1}$ with $r=\max\{1,m-1\}$. If $f\in W^{s-2,p}(\Omega)$ and $g\in W^{s-1/p,p}(\Gamma)$ with $s>1/p$, then $u\in W^{s,p}(\Omega)$ provided that $\Omega$ is of class $\mathscr{C}^{r,1}$ with $r=\max\{1,[s]\}$, where $[s]$ is the integer part of $s$. 

The 1980's saw many authors investing in the study of these problems when the domain $\Omega$ is only Lipschitz, where the situation is completely different (see for instance \cite{Cos, Dahl, Dahlberg, Fabes1, Fabes2, Jer1, Jer2, Verchota1}). Since then and to this day, these questions are still of great interest and many works are devoted to them (see for example \cite{Fabes3, GK, J-K, MMS, Mik, Mitrea1, Mitrea2, McL}).

Recall that when $\Omega $ is of class $\mathscr{C}^{1}$   and $1 < p < \infty$, for any $
f\in W^{-1,\, p}(\Omega)$ and for any $ g\in W^{1 - 1/p,\, p}(\Gamma),$ Problem $(\mathscr{L}_D)$ has a unique solution $
u\in W^{1,\, p}(\Omega).$ In the 80's, Ne\v{c}as posed the question of solving the problem $(\mathscr{L}^0_D)$ with the homogeneous boundary condition $g = 0$ on Lipschitz domains, when the RHS $f\in W^{-1,\, p}(\Omega)$. The answer to this question is given in Theorem A of the famous paper of Jerison and  Kenig \cite{J-K}: If $N \geq 3$, then for any $p > 3$, there is a Lipschitz domain $\Omega$ and $f\in \mathscr{C}^{\infty}(\overline{\Omega})$ such that the solution $u$ of Problem $(\mathscr{L}_D^0)$, that belongs to $H^1_0(\Omega)$,  does not belong to $W^{1,\, p}(\Omega)$ (for $N = 2$, the result is valid for any $p > 4$). However, for any Lipschitz bounded domain $\Omega$, there exists $q > 4$,  when $N = 2$ and  $q > 3$,  when $N \geq 2$, depending on $\Omega$,  such that if $q' < p < q,$ then the problem $(\mathscr{L}^0_D)$ has a unique solution $u\in W^{1,\, p}_0(\Omega)$ satisfying the estimate
$$
\Vert u \Vert_{W^{1,\, p}(\Omega)}\leq C \Vert f \Vert_{W^{-1,\, p}(\Omega)},
$$
%
(when $\Omega$ is $\mathscr{C}^1$, we can take $q = \infty$ as stated above). Observe that the exponents $4$ and $4/3$ in 2D,  respectively $3$ and $3/2$ in 3D,  corresponding respectively to the limit cases for the existence and the uniqueness of solutions in $W^{1,\, p}_0(\Omega)$, are conjugate. This is naturally due to the fact that the operator $\Delta : W^{1,\, p}_0(\Omega) \longrightarrow  W^{-1,\, p}(\Omega)$ is self-adjoint.

As an example of non-uniqueness, let us consider in the  2D case the following Lipschitz domain  for $1/2 < \alpha < 1$:
$$
\Omega = \{(r, \theta);\; 0 < r < 1,\quad 0 < \theta < \pi/\alpha\}.
$$
We can easily verify that the following function
\begin{equation}\label{exampleintro}
u(r, \theta) = (r^{-\alpha} - r^{\alpha})\mathrm{sin}(\alpha\theta)
\end{equation}
is harmonic in $\Omega$ with $u = 0$ on $\Gamma$ and $u\in W^{1, p}(\Omega)$ for any $p < 2/(\alpha + 1)$. Remark that when $\alpha$ is close to $1/2 $, then $\Omega$ is close to the unit disk and $2/(\alpha + 1)$ is close to  $4/3$. Note that the limit value 1/2 of $\alpha$ corresponds to the case of the cracked disk, which is not a Lipchitz domain.

In this work, we investigate the case where the data $f$ and $g$ are non smooth and the domain is only Lipschitz. These questions are not new and a lot of work has been done in this direction by many authors using various techniques since the 80's. If for regular domains almost every thing is elucidated, it is not the case for Lipschitz ones and for $s$ of the form
 $s = k + 1/2$, with $k$ integer and corresponding to limit cases. It is well known that this framework is delicate. So, we are particularly interested in investigating here the maximal regularity in these limit cases. In this direction, some new informations about the traces of functions belonging to $L^2(\Omega)$ (and even to $H^{-1}(\Omega)$) or to $H^{1/2}(\Omega)$ that satisfy an adequate additional property are very useful. One of the idea is to use the interpolation theory and the duality method.  Also, we prove new Ne\v{c}as' properties under the assumption of data satisfying less restrictive conditions than usual.

In the rest of this abstract,  $\Omega$ is a bounded Lipschitz domain in $\mathbb{R}^N$, with $N \geq 2$ and we give our main results. 

The first result is about traces of functions belonging to Sobolev spaces, which is crucial in the study of boundary value problems. We know that if $u\in H^s(\Omega)$ with $s > 1/2$ then the function $u$ has a trace which belongs to  $H^{s-1/2}(\Gamma)$. Moreover, if $v\in H^{1/2}(\Omega)$, in general this function $v$ may have no trace. In the following theorem, we see  that an additional condition on $\nabla v$  allows to obtain a trace for the function $v$.

\begin{theorem}[{\bf Trace Operator in $H^{1/2}(\Omega)$}]\label{TracesH1demigradH1demiprime1}  i) The linear mapping $\gamma_0: u \mapsto u_{\vert \Gamma}$ defined on $\mathscr{D}(\overline{\Omega})$ can be extended by continuity to a linear and continuous mapping, still denoted $\gamma_0$, from $E(\nabla;\, \Omega)$ into $L^2(\Gamma)$,  where
\begin{equation*}
\ds E(\nabla;\, \Omega)\ =\ \left\{\,v\in H^{1/2}(\Omega);\ \nabla v\in  [\textit{\textbf H}^{\, 1/2}(\Omega)]'\,\right\}.
\end{equation*}
ii) The Kernel of $\gamma_0$ is equal to $H^{1/2}_{00}(\Omega)$. In particular, we have the following property:
\begin{equation*}\label{traceH^{1/2}_{00}1}
v\in H^{1/2}_{00}(\Omega) \Longrightarrow v = 0 \quad in\; L^2(\Gamma). 
\end{equation*}
\end{theorem}
As a consequence we get immediately the following results: let 
$$
v\in H^{3/2}(\Omega)\quad \mathrm{with}\quad \nabla^2 v\in  [\textit{\textbf H}^{\, 1/2}(\Omega)]',
$$
then 
$$
v_{\vert\Gamma} \in H^1(\Gamma) \quad \mathrm{and}\quad \frac{\partial v}{\partial  \textit{\textbf n}}\in L^2(\Gamma).
$$
Moreover, we have the following property: 
\begin{equation*}\label{traceH^{3/2}_{00}1}
v\in H^{3/2}_{00}(\Omega) \Longrightarrow v = 0 \quad \mathrm{in}\; H^1(\Gamma)\quad  \mathrm{and}\quad \frac{\partial v}{\partial  \textit{\textbf n}}  = 0 \quad \mathrm{in}\; L^2(\Gamma).
\end{equation*}
The notation $\nabla^k v$ denotes the derivatives $\partial^\alpha v$ with $\vert \alpha \vert = k$. The vector fields and the spaces of vector fields are denoted by bold fonts.\medskip

The second result is about the $H^2$ regularity as well as the $H^3$ regularity for Dirichlet problem, which  is a question that has been addressed by many mathematicians. As above for solutions in $W^{1,\, p}_0(\Omega)$, we know that, for any $s > 3/2$,  there exist a Lipschitz domain $\Omega$ and $f \in \mathscr{C}^\infty(\overline{\Omega})$ such that the solution $u$ to the non homogeneous Dirichlet problem $(\mathscr{L}^0_D)$ does not belong to $H^s(\Omega)$ (see \cite{J-K} for example).  In the case of polygonal or polyedral domain $\Omega$, Grisvard in \cite{Gri1} and \cite{Gri2} gave for non negative integer $m$, a necessary and sufficient condition to obtain the solution in $H^{m+2}(\Omega)\cap H^1_0(\Omega)$ which can be written as follows: the RHS $f$ satisfies
\begin{equation*}
f\in H^m_0(\Omega) \quad \mathrm{and}\quad \langle v, f\rangle_{H^{-m}(\Omega) \times H^m_0(\Omega) }= 0\qquad \forall v\in \mathscr{H}^\circ_{H^{-m}(\Omega) } 
\end{equation*}
(see the notation below). In the following theorem, we extend this result to the general case of bounded Lipschitz domains. To start with, let us introduce the following notations. Let $X$ be a subspace of distributions in $\Omega$ and 
$$
\mathscr{H}_{X} = \{ v \in X ; \; \Delta v = 0 \; \mathrm{in}\,\,  \Omega  \},\quad \quad \mathscr{H}^\circ_{X} = \{ v \in \mathscr{H}_{X}  ; \;  v = 0\;  \mathrm{on}\; \Gamma \}
$$
and
$$
  (\mathscr{H}_{X})^\bot  = \{f\in X'; \; \forall \varphi\in \mathscr{H}_{X}, \;  \langle f,\, \varphi \rangle= 0\}, 
$$
with a similar definition for the orthogonal of $\mathscr{H}^\circ_{X} $. Observe that if $X = L^2(\Omega)$, then $X' = L^2(\Omega)$ and the duality brackets above can be replaced by an integral.


\begin{theorem}  [{\bf $H^2$ and $H^3$-Regularity}] \label{$H^2$-Regularity}   The operators
\begin{equation*}\label{DeltaIsoH2H10Ortho1}
\Delta : H^2(\Omega)\cap H^{1}_0(\Omega) \longrightarrow (\mathscr{H}^\circ_{L^{2}(\Omega)})^\bot
\end{equation*}
and
\begin{equation*}\label{DeltaIsoH3H10Ortho-1bis}
\Delta : H^3(\Omega)\cap H^{1}_0(\Omega) \longrightarrow  (\mathscr{H}^\circ_{H^{-1}(\Omega)})^\bot
\end{equation*}
are isomorphisms.
In particular we have the following inequalities: there exits a constant $C$ such that for any $v\in H^2(\Omega)\cap H^{1}_0(\Omega) $, we have
\begin{equation*}
\Vert v \Vert_{H^2(\Omega)} \leq C \Vert \Delta v \Vert_{L^2(\Omega)}.
\end{equation*}
and  for any $v\in H^3(\Omega)\cap H^{1}_0(\Omega)$ with $\Delta v\in H^{1}_0(\Omega) $, we have
\begin{equation*}
\Vert v \Vert_{H^3(\Omega)} \leq C \Vert \Delta v \Vert_{H^1(\Omega)}.
\end{equation*}
\end{theorem}
Observe that the space $\mathscr{H}^\circ_{L^{2}(\Omega)}$ (resp. $\mathscr{H}^\circ_{H^{-1}(\Omega)}$)  being the subspace of harmonic functions in $L^2(\Omega)$ (resp. in $H^{-1}(\Omega)$) vanishing on $\Gamma$, we need to give meaning to the traces of such functions or distributions. 

As an interesting consequence for the $H^2$ regularity of Problem $(\mathscr{L}_D^0)$, the operator
\begin{equation*}
\Delta : H^2(\Omega)\cap H^{1}_0(\Omega) \longrightarrow  L^{2}(\Omega)
\end{equation*}
is an isomorphism if and only if the kernel $\mathscr{H}^\circ_{L^{2}(\Omega)}$ is reduced to $0$.

\begin{remark} \upshape We know that this kernel  is reduced to $0$ in the case where $\Omega$ is of class $\mathscr{C}^{1, 1}$ or convex. One of the interesting questions that could be asked concerns the minimal regularity of the domain $\Omega$ to have a trivial kernel.
\end{remark}

Independently of the result above, we can show that the following operator
\begin{equation*}\label{DeltaIsoH20Ortho1}
\Delta : H^2_0(\Omega) \longrightarrow   (\mathscr{H}_{L^2(\Omega)})^\bot
\end{equation*}
is also an isomorphism. That means that for any $f\in L^2(\Omega) $ satisfying the compatibility condition 
\begin{equation*}\label{condorthoH1/21}
 \forall  \varphi \in \mathscr{H}_{L^{2}(\Omega)}, \quad \int_\Omega f \, \varphi   = 0,
\end{equation*} 
there exists a unique solution $u\in H^{2}_{0}(\Omega)$ such that $\Delta u = f$ in $\Omega$.  So, in addition to the boundary Dirichlet condition $u = 0$, the normal derivative of this solution also satisfies  $\frac{\partial u}{\partial  \textit{\textbf n}} = 0$. 

Another consequence of Theorem \ref{$H^2$-Regularity}  concerns the bi-Laplacian problem with the so-called Navier condition:
 \begin{equation}\label{BilapuDeltau=0b}
 \begin{cases}
 \ds\Delta^2 u\ =\ f \quad\textrm{ in }\ \Omega,\\
  \ds u\ = 0 \quad\textrm{ on }\ \Gamma,\\
 \Delta u \ = \ 0 \quad\textrm{ on }\ \Gamma.
 \end{cases}
 \end{equation}
We prove that for any
\begin{equation*}\label{corBilapuDeltau=0b}
f\in H^{-1}(\Omega)\quad \mathrm{such\, that} \quad (- \Delta)^{-1}f \in  (\mathscr{H}^\circ_{H^{-1}(\Omega)})^\bot,
\end{equation*}
there exists a unique function $u\in H^3(\Omega)$  solution to Problem \eqref{BilapuDeltau=0b}.\medskip

It is well known that if $f\in L^2(\Omega)$ (or even if $f\in H^{-s}(\Omega)$ for any $s < 1/2$), then there exists a unique solution $u\in H^{3/2}_0(\Omega)$ to Problem $(\mathscr{L}_D^0)$. But these assumptions on $f$ are too strong. So it would be interesting to characterize the range of  $H^{3/2}(\Omega)\cap H^1_0(\Omega)$ by the Laplacian operator. One of our main results, which was not proved yet so far as we know, is given by the next theorem:

%
\begin{theorem} [{\bf Solutions in $H^{3/2}_0(\Omega)$}] \label{IsoDeltaH3/2H1/21} i) The operators 
\begin{equation*}\label{DeltaIsoH3/2H1/21}
\Delta : H^{3/2}_0(\Omega) \longrightarrow [{H}^{1/2}_{00}(\Omega)]'\quad and \quad \Delta : H^{1/2}_{00}(\Omega) \longrightarrow  [H^{3/2}_0(\Omega)]'
\end{equation*}
are isomorphisms.\\
ii) For any $f\in [H^{1/2}(\Omega)]' $ satisfying the compatibility condition 
\begin{equation*}\label{condorthoH1/21}
 \forall  \varphi \in \mathscr{H}_{H^{1/2}(\Omega)}, \quad \langle f, \, \varphi  \rangle = 0,
\end{equation*} 
there exists a unique solution $u\in  H^{3/2}_{00}(\Omega)$ such that $\Delta u = f$ in $\Omega$.  In addition to the boundary Dirichlet condition $u = 0$, the normal derivative of this solution satisfies  $\frac{\partial u}{\partial  \textit{\textbf n}} = 0$. \\
\textit{\textbf {iii}})  By duality, for any $f\in \left[ H^{3/2}_{00}(\Omega)\right]'$, there exists $u\in H^{1/2}(\Omega)$, unique up to an element of $\mathscr{H}_{H^{1/2}(\Omega)}$, where 
$$
\mathscr{H}_{H^{1/2}(\Omega)}\ =\ \left\{ v\in H^{1/2}(\Omega); \ \ \Delta v=0  \;\; in\;\; \Omega\right\}.
$$
\end{theorem}

Using an interpolation argument, we deduce the following regularity results in $H^{s+1}(\Omega)$, with  $0 < s < 1$ and $s\not= 1/2$.

\begin{corollary}  [{\bf Solutions in $H^{s}(\Omega)$}]i) For any $0 < s < 1$, with $s\not= 1/2$, the operator
\begin{equation*}\label{DeltaIsoH3/2+sH1/2aa1b}
\Delta : H^{1+s}(\Omega)\cap H^{1}_0(\Omega)\longrightarrow  (\mathscr{H}^0_{H^{1-s}(\Omega)})^\bot
\end{equation*} 
is an isomorphism, where  $(\mathscr{H}^0_{H^{1-s}(\Omega)})^\bot$ denotes the orthogonal of the subspace of harmonic functions in $H^{s-1}(\Omega)$ vanishing on $\Gamma$.\\
ii) For any $1/2 < s < 1$,  the operator
\begin{equation*}\label{DeltaIsoH1+SH10s>1/2b}
\Delta : H^{1+s}_0(\Omega)\longrightarrow  (\mathscr{H}_{H^{1-s}(\Omega)})^\bot
\end{equation*} 
is an isomorphism
\end{corollary}

\begin{remark} \upshape Note that $\mathscr{H}^0_{H^{1-s}(\Omega)} = \{0\}$ if $0 < s < 1/2$ and then we have 
$$
H^{1+s}(\Omega) \cap H^1_0(\Omega) = H^{1+s}_0(\Omega) \quad \mathrm{and}\quad (\mathscr{H}^0_{H^{1-s}(\Omega)})^\bot  =  H^{s - 1}(\Omega).
$$
\end{remark}


%
%
%

We give now existence results in the case of boundary data in $L^2(\Gamma)$ or in $H^1(\Gamma)$. Using harmonic analysis techniques, many authors have established similar results (see \cite{Jer1} and \cite{J-K}). In the case where $g\in L^2(\Gamma)$, it is proved in \cite{Dahl77} the existence of a unique harmonic function such that $u$ tends  nontangentially to $g$ {\it a.e} on $\Gamma$ (see \cite{Dahl77})  and $u$ satisfies
$$
\Vert u^\star\Vert_{L^2(\Gamma)} \leq C \Vert g \Vert_{L^2(\Gamma)}.
$$
Here, the nontangential maximal function $u^\star$ is defined by
$$
\textit{\textbf z}\in \Gamma, \quad u^\star(\textit{\textbf z}) = \sup_{\textit{\textbf x}\, \in \, \Gamma(\textit{\textbf z})} \vert u(\textit{\textbf x})\vert,
$$
where $\Gamma(\textit{\textbf z})$ is a nontangential cone with vertex at $\textit{\textbf z}$, that is:
$$
\Gamma(\textit{\textbf z}) = \{\textit{\textbf x}\in \mathbb{R}^N; \; \vert \textit{\textbf x} - \textit{\textbf z}\vert < C\vert \varrho(\textit{\textbf x})\vert\},
$$
for a suitable constant $C > 1$. The notation $\varrho(\textit{\textbf x})$ denotes the distance from $\textit{\textbf x}\in \Omega$ to $\Gamma$. When $g\in H^1(\Gamma)$, there exists a unique harmonic function such that $u$ tends  nontangentially to $g$ {\it a.e} on $\Gamma$ (see \cite{Jer1})  and $u$ satisfies
$$
\Vert (\nabla u)^\star\Vert_{L^2(\Gamma)} \leq C \Vert g \Vert_{H^1(\Gamma)}.
$$
Our proofs, completely different, are essentially based on the first isomorphism given in Theorem \ref{IsoDeltaH3/2H1/21} and the following variant of  Ne\v{c}as' property: let 
\begin{equation*}\label{uH10DeltauH1/2prime1}
u\in H^1_0(\Omega)\quad \mathrm{ such \, that} \quad\Delta u  \in\ [H^{1/2}(\Omega)]',
\end{equation*}
then 
$$\
\ds\frac{\partial u}{\partial \textit{\textbf n}}\in L^2(\Gamma).
$$

\begin{theorem} [{\bf Homogeneous Problem in $H^{1/2}(\Omega)$ and in $H^{3/2}(\Omega)$}]\label{ThIsogL2Gamma1} i) For any $g\in L^{2}(\Gamma)$, Problem $(\mathscr{L}_D^H)$ has a unique solution $u\in H^{1/2}(\Omega)$. Moreover $\sqrt \varrho\, \nabla u  \in {\textit{\textbf L}}^{2}(\Omega)$ and the following estimate holds:
\begin{equation*}\label{d02-240118-e1c1}
\left|\left| u\right|\right|_{H^{1/2}(\Omega)} + \Vert\sqrt \varrho\, \nabla u \Vert_{\textit{\textbf L}^2(\Omega)}  \leq \ C\,\left|\left| g\right|\right|_{L^{2}(\Gamma)}.
\end{equation*}
The solution $u$ satisfies also the following property: for any positive integer $k$
$$
\varrho^{k + 1/2}\nabla^{k + 1}u \in \textit{\textbf L}^2(\Omega).
$$
ii) For any $g\in H^{1}(\Gamma)$, the problem $(\mathscr{L}_D^H)$ has a unique solution $u\in H^{3/2}(\Omega)$. Moreover $\sqrt \varrho \, \nabla^2 u \in {\textit{\textbf L}}^2(\Omega)$ and the following estimate holds:
\begin{equation*}\label{estimH3demiH1Gamma1}
\left|\left| u\right|\right|_{H^{3/2}(\Omega)} + \Vert\sqrt \varrho\, \nabla^2 u \Vert_{\textit{\textbf L}^2(\Omega)} \leq \ C\,\left|\left| g\right|\right|_{H^{1}(\Gamma)}.
\end{equation*}
The solution $u$ satisfies also the following property: for any positive integer $k$
$$
\varrho^{k + 1/2}\nabla^{k + 2}u \in \textit{\textbf L}^2(\Omega).
$$
\end{theorem}

\begin{remark}\upshape i) In addition to the existence and uniqueness result given in Point i) above, the sense of the boundary condition $u = g$ is as usual for boundary value problems the one given by the trace $L^2(\Gamma)$ and not only in the non-tangential sense as it has been found in the literature since the 80s.\\
ii) From Theorem \ref{IsoDeltaH3/2H1/21} and Point i) above, we deduce the following characterization: let $u\in H^{1/2}(\Omega)$, then
$$
u \in H^{1/2}_{00}(\Omega) \Longleftrightarrow \Delta u \in [H^{3/2}_0(\Omega)]' \quad \mathrm{and}\quad u = 0 \; \mathrm{on}\; \Gamma.
$$
With this characterization, we can also get the following: let $u\in H^{3/2}_0(\Omega)$, then
$$
u \in H^{3/2}_{00}(\Omega) \Longleftrightarrow \Delta u \in [H^{1/2}(\Omega)]' \quad \mathrm{and}\quad \frac{\partial u}{\partial \textit{\textbf n}}\ = 0 \; \mathrm{on}\; \Gamma.
$$
iii) In Point i) above, the properties  $\sqrt \varrho\, \nabla u  \in {\textit{\textbf L}}^{2}(\Omega)$ and \break$\varrho^{k + 1/2}\nabla^{k + 1}u \in \textit{\textbf L}^2(\Omega)$ are a direct consequence of the harmonicity of $u$.\\
iv) In a forthcoming paper, we will prove similar results in the case where the domain is of class $\mathscr{C}^{1,1}$ and the Dirichlet boundary condition $g$ belongs to $ H^2(\Gamma)$.
\end{remark}

According to Point i) of  Theorem \ref{ThIsogL2Gamma1},  it is natural to ask whether the reciprocal holds. Namely, if the properties   $ \sqrt \varrho\, \nabla u \in L^2(\Omega)$ and $u$ harmonic (so that  $u\in H^{1/2}(\Omega)$) imply that $u\in L^2(\Gamma)$? In \cite{Necas}, Ne\v{c}as proved the following property (see Theorem 2.2 Section 6): if $\varrho^{\alpha/p} u\in L^p(\Omega)$ and $\varrho^{\alpha/p}\nabla u \in L^p(\Omega)$, with $0 \leq \alpha < p - 1$, then $u_{\vert\Gamma} \in L^p(\Gamma)$ and 
\begin{equation}\label{inegNec}
\int_\Gamma \vert u \vert^p \leq C(\Omega) (\int_\Omega \varrho^\alpha\vert  u \vert^p + \int_\Omega \varrho^\alpha\vert \nabla u \vert^p)
\end{equation}
However, if $\alpha = p - 1$, the above inequality does not hold in general, as proved in a counter example with $\Omega= \,  ]0, 1/2[\,  \times\,  ]0, 1/2[ $. In particular, if $ \sqrt \varrho \, \nabla u \in L^2(\Omega)$, corresponding to the case $\alpha = 1$ and $p = 2$ and in which case we know that $u \in H^{1/2}(\Omega)$, the function $u$ in general has no trace. What about if in addition the function $u$ is harmonic?  What about if in addition the function $u$ is harmonic? In \cite{Dahlberg} (see also Corollary, Section 6 in \cite{Dahl}), using some specific properties of the distance to the boundary, the authors proved the following property: Let $u$ be a harmonic function in $\Omega$ that vanishes at some point $\textit{\textbf x}_0\in \Omega$, then
\begin{equation}\label{inegaltraceL2Gammabbis}
\int_\Gamma \vert u \vert^2 \leq C(\Omega)\int_\Omega \varrho\vert \nabla u \vert^2,
\end{equation}
where the constant $C(\Omega)$ depends only on the Lipschitz character of $\Omega.$
This result would then imply that the above inequality \eqref{inegNec} would be true in the critical case  $p = 2$ and $\alpha = 1$ when the function $u$ is in addition harmonic. 
Although the function $u$ is here harmonic, we will see below that the inequality \eqref{inegaltraceL2Gammabbis} is not true. Observe that the inequality \eqref{inegaltraceL2Gammabbis} implies the following property: let $u$ be a harmonic function in $\Omega$ satisfying $u(\textit{\textbf x}_0) = 0$ and $\nabla u(\textit{\textbf x}_0) = {\bf 0}$ at some point $\textit{\textbf x}_0\in \Omega$, then
\begin{equation}\label{inegaltraceL2Gammatera}
\Vert u \Vert_{H^1(\Gamma)} \leq C(\Omega)\Big(\int_\Omega \varrho\vert \nabla^2 u \vert^2\Big)^{1/2},
\end{equation}
where the constant $C(\Omega)$ depends only on the Lipschitz character of $\Omega.$ The following inequality would then be true for any harmonic function $u$:
\begin{equation}\label{inegaltraceL2Gammapent}
\Vert u \Vert_{H^1(\Gamma)} \leq C(\Omega)\Big(\int_\Omega \varrho\vert \nabla^2 u \vert^2 + \Vert u \Vert^2_{L^2(\Omega)}\Big)^{1/2}.
\end{equation}

The following proposition shows that Inequality  \eqref{inegaltraceL2Gammapent} is false and therefore Inequalities \eqref{inegaltraceL2Gammabbis}  and \eqref{inegaltraceL2Gammatera}   are also wrong.



\begin{e-proposition}  [{\bf Counter Example}]  \label{ConterexampleH1demitrace1}For any $\varepsilon > 0$, there is a Lipschitz domain $\Omega_\varepsilon \subset \mathbb{R}^2$ and a harmonic function $w_\varepsilon \in H^{3/2}(\Omega_\varepsilon)$ (with $\sqrt \varrho_\varepsilon\, \nabla^2 w_\varepsilon  \in L^2(\Omega_\varepsilon))$ such that the following family 
$$
(\Vert \varrho_\varepsilon \nabla^2 w_\varepsilon \Vert_{L^{2}(\Omega_\varepsilon)} + \Vert w_\varepsilon \Vert_{H^{3/2}(\Omega_\varepsilon)})_\varepsilon,
$$
is bounded with respect $\varepsilon$ and
\begin{equation*}
\Vert w_\varepsilon \Vert_{H^1(\Gamma_\varepsilon)} \rightarrow \infty \quad as \; \varepsilon\rightarrow  0.
\end{equation*}

\end{e-proposition}

Of course, by use Ne\v{c}as's Property below, we can replace the norm $\Vert w_\varepsilon \Vert_{H^1(\Gamma_\varepsilon)} $ by the norm $\Vert \nabla w_\varepsilon \Vert_{L^2(\Gamma_\varepsilon)} $ and then to get the following property: for any $\varepsilon > 0$, there is a Lipschitz domain $\Omega_\varepsilon \subset \mathbb{R}^2$ and a harmonic function $w_\varepsilon \in H^{1/2}(\Omega_\varepsilon)$ (with $\sqrt \varrho_\varepsilon\, \nabla w_\varepsilon  \in L^2(\Omega_\varepsilon))$ such that the following family 
$$
(\Vert \varrho_\varepsilon \nabla w_\varepsilon \Vert_{L^{2}(\Omega_\varepsilon)} + \Vert w_\varepsilon \Vert_{H^{1/2}(\Omega_\varepsilon)})_\varepsilon,
$$
is bounded with respect $\varepsilon$ and
\begin{equation*}
\Vert w_\varepsilon \Vert_{L^2(\Gamma_\varepsilon)} \rightarrow \infty \quad as \; \varepsilon\rightarrow  0.
\end{equation*}


We now give extensions of the classical Ne\v{c}as' property, that will be very useful for the study of the homogeneous Neumann problem and also for the  Dirichlet-to-Neumann operator for the Laplacian. Usually, $\Delta u$ is assumed to belong to $L^2(\Omega)$, which is a stronger condition than the one we take below.

\begin{theorem}  [{\bf Ne\v{c}as Property}] \label{d02-110118-th1a1}
Let 
$$
u\in H^1(\Omega)\quad \mathrm{ with} \quad \Delta u  \in\ [H^{1/2}(\Omega)]'.
$$
i) If $u\in H^1(\Gamma)$, then $\ds\frac{\partial u}{\partial \textit{\textbf n}}\in L^2(\Gamma)$ and we have the following estimate
\begin{equation*}\label{d02-100118-e2a1}
\ds\left|\left|\frac{\partial u}{\partial \textit{\textbf n}} \right|\right|_{ \textit{\textbf L}^2(\Gamma)}\ \leq\ C(\Omega)\left(\inf_{k\in \R}\left|\left |u + k \right |\right |_{H^1(\Gamma)} +\left|\left| \Delta u\right|\right|_{[H^{1/2}(\Omega)]'}\right), 
\end{equation*}
where the constant $C(\Omega)$ depends only on the Lipschitz character of $\Omega.$
\medskip

\noindent ii)  If $\ds\frac{\partial u}{\partial \textit{\textbf n}}\in L^2(\Gamma)$, then $u\in H^1(\Gamma)$ and we have the following estimate
\begin{equation*}\label{d02-100118-e3a1}
\ds \inf_{k\in \R}\left|\left|u + k \right|\right|_{H^1(\Gamma)}  \ \leq\ C(\Omega)\left( \left|\left|\frac{\partial u}{\partial \textit{\textbf n}} \right|\right|_{ \textit{\textbf L}^2(\Gamma)} +\left|\left|  \Delta u\right|\right|_{[H^{1/2}(\Omega)]')}\right),
\end{equation*}
where the constant $C(\Omega)$ depends only on the Lipschitz character of $\Omega.$\\
iii)  If $u\in H^1(\Gamma)$ or $\ds \frac{\partial u}{\partial\textit{\textbf n}}\in L^2(\Gamma)$, then $u\in H^{3/2}(\Omega)$.
\end{theorem}

\begin{remark}\upshape  In a forthcoming paper, using the properties of the  Dirichlet-to-Neumann operator, we will give different regularity results for the Neumann problem.
\end{remark}

In this work we investigate also the case of solutions in the following weighted Sobolev spaces: 
 $$
 \mathscr{T}^2_{-3/2}(\Omega) = \{v \in H^{3/2}_0(\Omega); \; \sqrt \varrho \, \nabla^2 v \in L^2(\Omega)\}
 $$
 and 
 $$
 L^2_{\sqrt \varrho}(\Omega) = \{v\in \mathscr{D}'(\Omega); \; \sqrt \varrho \,  v \in L^2(\Omega) \}
 $$
 and note that $\mathscr{D}(\Omega)$ is dense in each of them.  So that their respective duals denoted by $\mathscr{T}^{-2}_{3/2}(\Omega)$ and $L^2_{1/\sqrt\varrho}(\Omega)$ are subspaces of distributions.

%
\begin{theorem} [{\bf Solutions in Weighted Sobolev Spaces}] \label{theoisospeciaux}The following operators, each being dual to the other, 
\begin{equation*}\label{isoDeltaF3/2etEthetaMR}
\Delta:\ \ \ \mathscr{T}^2_{-3/2}(\Omega)\longrightarrow \ \ \ L^2_{\sqrt \varrho}(\Omega) \quad and \quad \Delta : L^2_{1/\sqrt\varrho}(\Omega)\longrightarrow \mathscr{T}^{-2}_{3/2}(\Omega)
\end{equation*}
are isomorphisms.   In particular, we have the following property:
\begin{equation}\label{uniquenesscriteriumMR}
v\in L^2_{1/\sqrt\varrho}(\Omega) \quad and \quad \Delta v\in  [H^{3/2}_0(\Omega)]'\quad \Longrightarrow \quad v \in H^{1/2}_{00}(\Omega).
\end{equation}
\end{theorem}

\begin{remark}\upshape  i) The first isomorphism above can be considered as a regularity result with respect to the first isomorphism of Theorem \ref{IsoDeltaH3/2H1/21}, Point i) since we have the following embeddings:
$$
\mathscr{T}^2_{-3/2}(\Omega)\hookrightarrow H^{3/2}_0(\Omega) \quad \mathrm{and} \quad L^2_{\sqrt \varrho}(\Omega)\hookrightarrow [H^{1/2}_{00}(\Omega)]'. 
$$
ii) Similarly the second isomorphism can be considered as an extension of the second isomorphism given by Theorem \ref{IsoDeltaH3/2H1/21}.\\
iii) In the book of Lions-Magenes, Chapter 2 \cite{Lions}, where the domain $\Omega$ is of class $\mathscr{C}^\infty$, (see also Ne\v{c}as Chapter 6 \cite{Necas}), the concept of very weak solutions in $H^{-s}(\Omega)$ for Problem  $(\mathscr{L}_D^0)$, with $s \geq 0$, has been widely developed. In particular for $s = 0$, they proved that for any $f\in \ds[\Xi^2(\Omega)]'$, where
$$
\ds\Xi^2(\Omega)\ =\ \left\{ v\in L^2(\Omega);\ \ \varrho^{|\alpha|}D^\alpha v\in \textbf{\textit L}^2(\Omega),\ \ |\alpha|\leq 2 \right\},
$$
there exists a unique solution $u\in L^2(\Omega)$ to Problem  $(\mathscr{L}_D^0)$. The second isomorphism above, where $\Omega$ is  only Lipschitz, improves significantly this result since
$$
\ds[\Xi^2(\Omega)]' \hookrightarrow \mathscr{T}^{-2}_{3/2}(\Omega)\quad \mathrm{and} \quad L^2_{1/\sqrt\varrho}(\Omega)\hookrightarrow L^2(\Omega). 
$$
iv)  In other words for \eqref{uniquenesscriteriumMR}, if $v$ is an harmonic function non identically equal to zero and vanishing on the boundary, then it cannot belong to the space $L^2_{1/\sqrt\varrho}(\Omega)$ as we can see in the above example \eqref{exampleintro}. Indeed, we easily verify that if $u$ belongs to $L^2_{1/\varrho^{\beta}}(\Omega)$, then necessarily $\beta < 1/2$ since $2\alpha > 1$.
 \end{remark}

{\small{\textit{Acknowledgments.} 
We are deeply indebted towards Professor  Martin Costabel and Professor David Jerison.
Our exchanges have been extremely  fruitful for us.
Thanks a lot for them.
}

\end{document}